\documentclass[11pt]{article}
\usepackage{amsfonts}
\usepackage{amsfonts}
\usepackage{amsfonts}
\usepackage{mathrsfs}
\usepackage{bm}
\usepackage{amssymb,amsmath} % font used for R in Real numbers
 \allowdisplaybreaks

\catcode`!=11
\let\!int\int \def\int{\displaystyle\!int}
\let\!lim\lim \def\lim{\displaystyle\!lim}
\let\!sum\sum \def\sum{\displaystyle\!sum}
\let\!sup\sup \def\sup{\displaystyle\!sup}
\let\!inf\inf \def\inf{\displaystyle\!inf}
\let\!cap\cap \def\cap{\displaystyle\!cap}
\let\!max\max \def\max{\displaystyle\!max}
\let\!min\min \def\min{\displaystyle\!min}
\let\!frac\frac \def\frac{\displaystyle\!frac}
\catcode`!=12

\let\oldsection\section
\renewcommand\section{\setcounter{equation}{0}\oldsection}

\allowdisplaybreaks
\def\pf{\it{Proof.}\rm\quad}

\def\N{\mathbb{N}}

\newtheorem{thm}{Theorem}[section]
\newtheorem{lem}[thm]{Lemma}
\newtheorem{cor}[thm]{Corollary}

\begin{document}
%%%%%%%%%%%%%%%%%%%% title %%%%%%%%%%%%%%%%%%%%%%%%%%%%%%%%%%%%%%%%%%%%%%%%
\title {\bf   Identities for the shifted harmonic numbers and binomial coefficients}
\author{
{Ce Xu\thanks{Email: 15959259051@163.com (C. Xu)}}\\[1mm]
\small School of Mathematical Sciences, Xiamen University\\
\small Xiamen
361005, P.R. China}

\date{}
\maketitle \noindent{\bf Abstract }  We develop new closed form representations of sums of $(n+\alpha)$th shifted harmonic numbers and reciprocal binomial coefficients through $\alpha$th shifted harmonic numbers and Riemann zeta function with positive integer arguments. Some interesting new consequences and illustrative examples are considered.
\\[2mm]
\\[2mm]

\section{Introduction}
Let $\N:=\{1,2,3,\ldots\}$ be the set of natural numbers, and $\N_0:=\N\cup\{0\}$. In this paper we will develop identities, closed form representations of shifted harmonic numbers and reciprocal binomial coefficients of the form:
\[{W^{(l)}_{k,r}}\left( {p,m,\alpha } \right): = \sum\limits_{n = 1}^\infty  {\frac{({H_{n + \alpha }^{(m)})^l}}{{{n^p}\left( {\begin{array}{*{20}{c}}
   {n + k + r}  \\
   k  \\
\end{array}} \right)}}},\ (r,k\in \N_0,\ \alpha \notin  \N^-:=\{-1,-2,\ldots\}),\tag{1.1}\]
for $m,l\in \{1,2,3\}, l+m\in\{2,3,4\}, p\in \{0,1\}$ with $p+k\geq 2$, where $H^{(m)}_\alpha$ stands for the $\alpha$-th generalized shifted harmonic number defined by
\[H_\alpha ^{\left( m \right)}: = \frac{{{{\left( { - 1} \right)}^{m - 1}}}}{{\left( {m - 1} \right)!}}\left( {{\psi ^{\left( {m - 1} \right)}}\left( {\alpha  + 1} \right) - {\psi ^{\left( {m - 1} \right)}}\left( 1 \right)} \right),\;2 \le m \in \N,\tag{1.2}\]
\[{H_\alpha } = H_\alpha ^{\left( 1 \right)}: = \psi \left( {\alpha  + 1} \right) + \gamma ,\tag{1.3}\]
here $\psi \left( z \right)$ is digamma function (or called Psi function ) defined by
\[\psi \left( z \right) := \frac{d}{{dz}}\left( {\ln \Gamma \left( z \right)} \right) = \frac{{\Gamma '\left( z \right)}}{{\Gamma \left( z \right)}}\]
and $\psi \left( z \right)$ satisfy the following relations in the forms
\[\psi \left( z \right) =  - \gamma  + \sum\limits_{n = 0}^\infty  {\left( {\frac{1}{{n + 1}} - \frac{1}{{n + z}}} \right)} ,\;z\notin  \N^-_0:=\{0,-1,-2\ldots\}, \]
\[{\psi ^{\left( n \right)}}\left( z \right) = {\left( { - 1} \right)^{n + 1}}n!\sum\limits_{k = 0}^\infty  {1/{{\left( {z + k} \right)}^{n + 1}}}, n\in \N,\]
\[\psi \left( {x + n} \right) = \frac{1}{x} + \frac{1}{{x + 1}} +  \cdots  + \frac{1}{{x + n - 1}} + \psi \left( x \right),\;n = 1,2,3, \ldots .\]
From the definition of Riemann zeta function and Hurwitz zeta function, we know that
\[{\psi ^{\left( n \right)}}\left( 1 \right) = {\left( { - 1} \right)^{n + 1}}n!\zeta \left( {n + 1} \right),{\psi ^{\left( n \right)}}\left( z \right) = {\left( { - 1} \right)^{n + 1}}n!\zeta \left( {n + 1,z} \right).\]
The Riemann zeta function and Hurwitz zeta function are defined by
\[\zeta(s):=\sum\limits_{n = 1}^\infty {\frac {1}{n^{s}}},\Re(s)>1,\tag{1.4}\]
and
\[\zeta \left( {s,\alpha  + 1} \right): = \sum\limits_{n = 1}^\infty  {\frac{1}{{{{\left( {n + \alpha } \right)}^s}}}} ,\left( {{\mathop{\Re}\nolimits} \left( s \right) > 1,\alpha  \notin \N^ - } \right).\tag{1.5}\]
Therefore, in the case of non integer values we may write the generalized shifted harmonic numbers
in terms of zeta functions
\[H_\alpha ^{\left( m \right)} = \zeta \left( m \right) - \zeta \left( {m,\alpha  + 1} \right),\;\alpha \notin  \N^-,\ 2\leq m\in \N,\tag{1.6}\]
and for $m=1$,
\[{H_\alpha } \equiv H_\alpha ^{\left( 1 \right)} = \sum\limits_{k = 1}^\infty  {\left( {\frac{1}{k} - \frac{1}{{k + \alpha }}} \right)} .\tag{1.7}\]
Here $\Gamma \left( z \right) := \int\limits_0^\infty  {{e^{ - t}}{t^{z - 1}}dt} ,\;{\mathop{\Re}\nolimits} \left( z \right) > 0$ is called gamma function, and $\gamma$ denotes the Euler-Mascheroni constant defined by
\[\gamma  := \mathop {\lim }\limits_{n \to \infty } \left( {\sum\limits_{k = 1}^n {\frac{1}{k}}  - \ln n} \right) =  - \psi \left( 1 \right) \approx {\rm{ 0 }}{\rm{. 577215664901532860606512 }}....\]
The evaluation of the polygamma function $\psi^{(n)}(\frac{p}{q})$ at rational values of the argument
can be explicitly done via a formula as given by K$\ddot{o}$lbig [14], or Choi and Cvijovi$\acute{c}$ [10] in
terms of the Polylogarithmic or other special functions. Some specific values are listed
in the books [1,18,25]. For example, George E. Andrews, Richard Askey and Ranjan Roy [1], or H.M. Srivastava and J. Choi [25] gave the following Gauss's formula
\begin{align*}
\psi \left( {\frac{p}{q}} \right) =& 2\sum\limits_{k = 1}^{\left[ {(q - 1)/2} \right]} {\cos \left( {\frac{{2kp\pi }}{q}} \right)\ln \left( {2\sin \frac{{k\pi }}{q}} \right)}  + {r_q}\left( p \right)\\&  - \gamma  - \frac{\pi }{2}\cot \frac{{p\pi }}{q} - \ln q ,\;\ \ \ \ \ \ \ \ \ \left( {0 < p < q;p,q \in \N} \right),\tag{1.8}
\end{align*}
where $[x]$ denotes the greatest integer $\leq x$, and if $q$ is even, ${r_q}\left( p \right) = {\left( { - 1}\right)^p}\ln 2$; if $q$ is odd, ${r_q}\left( p \right) = 0$. For details and historical introductions, please see [1, 3, 4, 10, 14, 18, 25] and references therein.
From (1.2), (1.3) and (1.8), we can obtain some specific values of shifted harmonic numbers:
\[{H_{1/2}} = 2 - 2\ln 2,{H_{3/2}} = \frac{8}{3} - 2\ln 2,H_{1/2}^{\left( 2 \right)} = 4 - 2\zeta \left( 2 \right),H_{3/2}^{\left( 2 \right)} = \frac{{40}}{9} - 2\zeta \left( 2 \right),{H_{5/2}} = \frac{{46}}{{15}} - 2\ln 2.\]
Letting $\alpha$ approach $n$ ($n$ is a positive integer) in (1.6) and (1.7), then the shifted harmonic numbers are  reducible to classical harmonic number defined by
$${H_n} \equiv H_n^{\left( 1 \right)}:=\sum\limits_{j=1}^n\frac {1}{j},\ H^{(m)}_n:=\sum\limits_{j=1}^n\frac {1}{j^m},\ n, m \in \N.\eqno(1.9)$$
While there are many results for sums of harmonic numbers (or shifted harmonic numbers) with positive terms, for
example we know that [20,23]
\[\sum\limits_{n = 1}^\infty  {\frac{{H_n^{ 2}}}{{\left( {\begin{array}{*{20}{c}}
   {n + k}  \\
   k  \\
\end{array}} \right)}}}  = \frac{k}{{k - 1}}\left( {\zeta \left( 2 \right) - H_{k - 1}^{\left( 2 \right)} + \frac{2}{{{{\left( {k - 1} \right)}^2}}}} \right),\;2 \le k \in \N,\]
\[\sum\limits_{n = 1}^\infty  {\frac{{H_{n + 1/4}^{\left( 2 \right)}}}{{n\left( {\begin{array}{*{20}{c}}
   {n + 4}  \\
   4  \\
\end{array}} \right)}}}  = 152 - \frac{{140}}{3}G - \frac{{128}}{3}\ln 2 - \frac{{64}}{9}\pi  - \frac{{69}}{2}\zeta \left( 2 \right),\]
where $G$ is Catalan¡¯s constant, defined by
\[G = \sum\limits_{n = 1}^\infty  {\frac{{{{\left( { - 1} \right)}^{n - 1}}}}{{{{\left( {2n - 1} \right)}^2}}}}  \approx {\rm{ 0}}{\rm{.915965}}....\]
Further work in the summation of harmonic numbers and binomial coefficients has also been done by Sofo [18-24].
In this paper, we will prove that the series of (1.1) can be expressed as a rational linear combination of products of zeta values and shifted harmonic numbers. Next, we give two lemmas.
The following lemma will be useful in the development of the main theorems.
\begin{lem} For integers $k,r\geq 0$ and $\alpha\not\in \N^-$, then have
\[\sum\limits_{n = 1}^\infty  {\frac{{f\left( {n,\alpha } \right)}}{{\left( {n + r} \right)\left( {n + k} \right)}}}  = \frac{{k - \alpha }}{{k - r}}\sum\limits_{n = 1}^\infty  {\frac{{f\left( {n,\alpha } \right)}}{{\left( {n + k} \right)\left( {n + \alpha } \right)}}}  + \frac{{\alpha  - r}}{{k - r}}\sum\limits_{n = 1}^\infty  {\frac{{f\left( {n,\alpha } \right)}}{{\left( {n + r} \right)\left( {n + \alpha } \right)}}} ,\tag{1.10}\]
where the function $f\left( {n,\alpha } \right)$ satisfy the following relation
\[\mathop {\lim }\limits_{n \to \infty } {n^\beta }f\left( {n,\alpha } \right) = c,\;\beta  >  - 1,\]
$c$ is a constant.
\end{lem}
\pf This lemma is almost obvious.\hfill$\square$
\begin{lem} For integer $p>0$ and $\alpha \neq 0,-1,-2,\ldots$, we have
\[\int\limits_0^1 {{x^{\alpha  - 1}}{\rm{L}}{{\rm{i}}_p}\left( x \right)dx}  = \sum\limits_{i = 1}^{p - 1} {\frac{{{{\left( { - 1} \right)}^{i-1}}}}{{{\alpha ^i}}}\zeta \left( {p + 1 - i} \right)}  - {\left( { - 1} \right)^p}\frac{{{H_\alpha }}}{{{\alpha ^p}}},\tag{1.11}\]
where ${\rm Li}{_p}\left( x \right)$ is polylogarithm function defined for $\left| x \right| \le 1$ by
\[{\rm Li}{_p}\left( x \right) = \sum\limits_{n = 1}^\infty  {\frac{{{x^n}}}{{{n^p}}}}, \Re(p)>1 .\tag{1.12}\]
\end{lem}
\pf It is obvious that we can rewrite the integral on the left hand side of (1.11) as
\[\int\limits_0^1 {{x^{\alpha  - 1}}{\rm{L}}{{\rm{i}}_p}\left( x \right)dx}  = \sum\limits_{n = 1}^\infty  {\frac{1}{{{n^p}\left( {n + \alpha } \right)}}}.\tag{1.13} \]
By using the partial fraction decomposition
\[\frac{1}{{{n^p}\left( {n + \alpha } \right)}} = \sum\limits_{i = 1}^{p - 1} {\frac{{{{\left( { - 1} \right)}^{i-1}}}}{{{\alpha ^i}}} \cdot \frac{1}{{{n^{p + 1 - i}}}}}  - {\left( { - 1} \right)^p}\frac{1}{{{\alpha ^p}}}\left( {\frac{1}{n} - \frac{1}{{n + \alpha }}} \right),\tag{1.14}\]
and combining (1.7) with (1.13), we deduce the desired result. This completes the proof
of Lemma 1.2.\hfill$\square$
\section{Shifted harmonic number identities}
In this section, we will establish some explicit relationships that involve shifted harmonic numbers and the following type sums
\[\sum\limits_{n = 1}^\infty  {\frac{{H_{n + \alpha }^{\left( m \right)}}}{{\left( {n + r} \right)\left( {n + k} \right)}}} ,\;\left( {r,k\in \N_0,\ r \neq k,\alpha  \notin \N^-,m\in \{1,2\}} \right).\]
We now prove the following theorems.
\begin{thm} For integers $m,k\geq 1$ and $\Re(\alpha)>0$, then we have the following recurrence relation
\begin{align*}
&I\left( {\alpha ,m,k} \right) = \sum\limits_{i = 0}^{m - 1} {\left( {\begin{array}{*{20}{c}}
   {m - 1}  \\
   i  \\
\end{array}} \right)\left( {m - i - 1} \right)!\frac{{{{\left( { - 1} \right)}^{m - i}}}}{{{\alpha ^{m - i}}}}I\left( {\alpha ,i,k} \right)} \\
&\quad\  + \sum\limits_{j = 0}^{k - 1} {\left( {\begin{array}{*{20}{c}}
   k  \\
   j  \\
\end{array}} \right){{\left( { - 1} \right)}^{m + k - j}}\left( {m + k - j - 1} \right)!H_\alpha ^{\left( {m + k - j} \right)}I\left( {\alpha ,0,j} \right)} \\
&\quad\  - \sum\limits_{j = 0}^{k - 1} {\left( {\begin{array}{*{20}{c}}
   k  \\
   j  \\
\end{array}} \right){{\left( { - 1} \right)}^{m + k - j}}\left( {m + k - j - 1} \right)!\zeta \left( {m + k - j} \right)I\left( {\alpha ,0,j} \right)} \\
&\quad\  + \sum\limits_{i = 1}^{m - 1} {\sum\limits_{j = 0}^{k - 1} {\left( {\begin{array}{*{20}{c}}
   {m - 1}  \\
   i  \\
\end{array}} \right)\left( {\begin{array}{*{20}{c}}
   k  \\
   j  \\
\end{array}} \right){{\left( { - 1} \right)}^{m + k - i - j}}\left( {m + k - i - j - 1} \right)!H_\alpha ^{\left( {m + k - i - j} \right)}I\left( {\alpha ,i,j} \right)} } \\
&\quad\  - \sum\limits_{i = 1}^{m - 1} {\sum\limits_{j = 0}^{k - 1} {\left( {\begin{array}{*{20}{c}}
   {m - 1}  \\
   i  \\
\end{array}} \right)\left( {\begin{array}{*{20}{c}}
   k  \\
   j  \\
\end{array}} \right){{\left( { - 1} \right)}^{m + k - i - j}}\left( {m + k - i - j - 1} \right)!\zeta \left( {m + k - i - j} \right)I\left( {\alpha ,i,j} \right)} }.\tag{2.1}
\end{align*}
where $I\left( {\alpha ,m,k} \right)$ is defined by the integral
\[I\left( {\alpha ,m,k} \right): = \int\limits_0^1 {{x^{\alpha  - 1}}{{\ln }^m}x{{\ln }^k}\left( {1 - x} \right)} dx.\tag{2.2}\]
with
\[I\left( {\alpha ,0,0} \right) = \frac{1}{\alpha },I\left( {\alpha ,i,0} \right) = {\left( { - 1} \right)^i}i!\frac{1}{{{\alpha ^{i + 1}}}}.\]
\end{thm}
\pf Applying the definition of Beta function ${B\left( {\alpha ,\beta } \right)}$, we can find that
\[I\left( {\alpha ,m,k} \right): = \int\limits_0^1 {{x^{\alpha  - 1}}{{\ln }^m}x{{\ln }^k}\left( {1 - x} \right)} dx = {\left. {\frac{{{\partial ^{m + k}}B\left( {\alpha ,\beta } \right)}}{{\partial {\alpha ^m}\partial {\beta ^k}}}} \right|_{\beta  = 1}},\tag{2.3}\]
where the Beta function is defined by
\[B\left( {\alpha,\beta} \right) := \int\limits_0^1 {{x^{\alpha - 1}}{{\left( {1 - x} \right)}^{\beta - 1}}dx}  = \frac{{\Gamma \left( \alpha \right)\Gamma \left( \beta \right)}}{{\Gamma \left( {\alpha + \beta} \right)}},\;{\mathop{\Re}\nolimits} \left( \alpha \right) > 0,{\mathop{\Re}\nolimits} \left( \beta \right) > 0.\tag{2.4}\]
By using (2.4) and the definition of $\psi (x)$, it is obvious that
\[\frac{{\partial B\left( {\alpha ,\beta } \right)}}{{\partial \alpha }} = B\left( {\alpha ,\beta } \right)\left[ {\psi \left( \alpha  \right) - \psi \left( {\alpha  + \beta } \right)} \right].\]
Therefore, differentiating $m-1$ times this equality, we can deduce that
\[\frac{{{\partial ^m}B\left( {\alpha ,\beta } \right)}}{{\partial {\alpha ^m}}} = \sum\limits_{i = 0}^{m - 1} {\left( {\begin{array}{*{20}{c}}
   {m - 1}  \\
   i  \\
\end{array}} \right)\frac{{{\partial ^i}B\left( {\alpha ,\beta } \right)}}{{\partial {\alpha ^i}}}} \cdot\left[ {{\psi ^{\left( {m - i - 1} \right)}}\left( \alpha  \right) - {\psi ^{\left( {m - i - 1} \right)}}\left( {\alpha  + \beta } \right)} \right].\tag{2.5}\]
Since $B(\alpha,\beta)=B(\beta,\alpha)$, then we also have
\[\frac{{{\partial ^m}B\left( {\alpha ,\beta } \right)}}{{\partial {\beta ^m}}} = \sum\limits_{i = 0}^{m - 1} {\left( {\begin{array}{*{20}{c}}
   {m - 1}  \\
   i  \\
\end{array}} \right)\frac{{{\partial ^i}B\left( {\alpha ,\beta } \right)}}{{\partial {\beta ^i}}}} \cdot\left[ {{\psi ^{\left( {m - i - 1} \right)}}\left( \beta  \right) - {\psi ^{\left( {m - i - 1} \right)}}\left( {\beta  + \alpha } \right)} \right].\tag{2.6}\]
Putting $\beta=1$ in (2.6) and combining (1.2) (2.3), we arrive at the conclusion that
\[I\left( {\alpha ,0,m} \right){\rm{ = }}\sum\limits_{i = 0}^{m - 1} {{{\left( { - 1} \right)}^{m - i}}\left( {m - i - 1} \right)!\left( {\begin{array}{*{20}{c}}
   {m - 1}  \\
   i  \\
\end{array}} \right)I\left( {\alpha ,0,i} \right)} H_\alpha ^{\left( {m - i} \right)}.\tag{2.7}\]
Furthermore, by using (2.5), the following identity is easily derived
\begin{align*}
\frac{{{\partial ^{m + k}}B\left( {\alpha,\beta} \right)}}{{\partial {\alpha^m}\partial {\beta^k}}}
 &=\frac{{{\partial ^k}}}{{\partial {\beta^k}}}\left( {\frac{{{\partial ^m}B\left( {\alpha,\beta} \right)}}{{\partial {\alpha^m}}}} \right)
\nonumber \\
           &=\sum\limits_{i = 0}^{m - 1} {\left( {\begin{array}{*{20}{c}}
   {m - 1}  \\
   i  \\
\end{array}} \right)\frac{{{\partial ^{i + k}}B\left( {\alpha,\beta} \right)}}{{\partial {\alpha^i}\partial {\beta^k}}} \cdot } \left[ {{\psi ^{\left( {m - i - 1} \right)}}\left(\alpha \right) - {\psi ^{\left( {m - i - 1} \right)}}\left( {\alpha + \beta} \right)} \right]
\nonumber \\
           &\quad \  - \sum\limits_{j = 0}^{k - 1} {\left( {\begin{array}{*{20}{c}}
   k  \\
   j  \\
\end{array}} \right)} \frac{{{\partial ^{ j}}B\left( {\alpha,\beta} \right)}}{{\partial {\beta^j}}}{\psi ^{\left( {m + k - j - 1} \right)}}\left( {\alpha + \beta} \right)
\nonumber \\
           &\quad \  - \sum\limits_{i = 1}^{m - 1} {\sum\limits_{j = 0}^{k - 1} {\left( {\begin{array}{*{20}{c}}
   {m - 1}  \\
   i  \\
\end{array}} \right)\left( {\begin{array}{*{20}{c}}
   k  \\
   j  \\
\end{array}} \right)} \frac{{{\partial ^{i + j}}B\left( {\alpha,\beta} \right)}}{{\partial {\alpha^i}\partial {\beta^j}}}} {\psi ^{\left( {m + k - i - j - 1} \right)}}\left( {\alpha + \beta} \right)
.\tag{2.8}
\end{align*}
From (1.2) and (1.5), we know that if $\beta=1$, then we have
\[{\psi ^{\left( {m - i - 1} \right)}}\left( \alpha  \right) - {\psi ^{\left( {m - i - 1} \right)}}\left( {\alpha  + 1} \right) = {\left( { - 1} \right)^{m - i}}\left( {m - i - 1} \right)!\frac{1}{{{\alpha ^{m - i}}}},\tag{2.9}\]
\[{\psi ^{\left( {m + k - j - 1} \right)}}\left( {\alpha  + 1} \right) = {\left( { - 1} \right)^{m + k - j}}\left( {m + k - j - 1} \right)!\left( {\zeta \left( {m + k - j} \right) - H_\alpha ^{\left( {m + k - j} \right)}} \right).\tag{2.10}\]
Hence, taking $\beta=1$ in (2.8), then substituting (2.9) and (2.10) into (2.8) respectively, we can obtain (2.1). The proof of Theorem 2.1 is finished.\hfill$\square$\\
From (2.1) and (2.7), we can get the following identities: for $\alpha>0$,
\begin{align*}
&I\left( {\alpha ,0,1} \right) = \int\limits_0^1 {{x^{\alpha  - 1}}\ln \left( {1 - x} \right)dx}  =  - \frac{{{H_\alpha }}}{\alpha },\tag{2.11}\\
&I\left( {\alpha ,0,2} \right) = \int\limits_0^1 {{x^{\alpha  - 1}}{{\ln }^2}\left( {1 - x} \right)dx}  = \frac{{H_\alpha ^2 + H_\alpha ^{\left( 2 \right)}}}{\alpha },\tag{2.12}\\
&I\left( {\alpha ,1,1} \right) = \int\limits_0^1 {{x^{\alpha  - 1}}\ln x\ln \left( {1 - x} \right)dx}  = \frac{{{H_\alpha }}}{{{\alpha ^2}}} - \frac{{\zeta \left( 2 \right) - H_\alpha ^{\left( 2 \right)}}}{\alpha },\tag{2.13}\\
&I\left( {\alpha ,0,3} \right) = \int\limits_0^1 {{x^{\alpha  - 1}}{{\ln }^3}\left( {1 - x} \right)dx}  =  - \frac{{H_\alpha ^3 + 3{H_\alpha}H_\alpha ^{\left( 2 \right)} + 2H_\alpha ^{\left( 3 \right)}}}{\alpha },\tag{2.14}\\
&I\left( {\alpha ,1,2} \right) = \int\limits_0^1 {{x^{\alpha  - 1}}\ln x{{\ln }^2}\left( {1 - x} \right)dx}  =  - \frac{{H_\alpha ^2 + H_\alpha ^{\left( 2 \right)}}}{{{\alpha ^2}}} + 2\frac{{\zeta \left( 3 \right) - H_\alpha ^{\left( 3 \right)}}}{\alpha } + 2\frac{{\zeta \left( 2 \right) - H_\alpha ^{\left( 2 \right)}}}{\alpha }{H_\alpha }.\tag{2.15}
\end{align*}
\begin{thm} For integers $r,k\in \N_0 \ (r\neq k)$ and $\alpha> max\{k,r\}$, we have
\begin{align*}
\sum\limits_{n = 1}^\infty  {\frac{{{H_{n + \alpha}}}}{{\left( {n + r} \right)\left( {n + k} \right)}}}  =& \frac{{k - \alpha }}{{k - r}}\left\{ {\frac{{H_{\alpha  - k}^2 + H_{\alpha  - k}^{\left( 2 \right)}}}{{\alpha  - k}} - \sum\limits_{j = 1}^k {\frac{{{H_{\alpha  + j - k}}}}{{j\left( {\alpha  + j - k} \right)}}} } \right\}\\
 &+ \frac{{\alpha  - r}}{{k - r}}\left\{ {\frac{{H_{\alpha  - r}^2 + H_{\alpha  - r}^{\left( 2 \right)}}}{{\alpha  - r}} - \sum\limits_{j = 1}^r {\frac{{{H_{\alpha  + j - r}}}}{{j\left( {\alpha  + j - r} \right)}}} } \right\}.\tag{2.16}
\end{align*}
\end{thm}
\pf Replacing $\alpha$ by $n+\alpha$ in (2.11), then multiplying it by $(n+k)^{-1}$ and summing with respect to $n$,  we conclude that
\begin{align*}
\sum\limits_{n = 1}^\infty  {\frac{{{H_{n + \alpha }}}}{{\left( {n + k} \right)\left( {n + \alpha } \right)}}}  =&  - \sum\limits_{n = 1}^\infty  {\frac{1}{{n + k}}\int\limits_0^1 {{x^{n + \alpha  - 1}}\ln \left( {1 - x} \right)} } dx\\
 = &\sum\limits_{j = 1}^k {\frac{1}{j}\int\limits_0^1 {{x^{\alpha  + j - k - 1}}\ln \left( {1 - x} \right)dx} }  + \int\limits_0^1 {{x^{\alpha  - k - 1}}{{\ln }^2}\left( {1 - x} \right)dx}.\tag{2.17}
\end{align*}
The relations (2.11), (2.12) and (2.17) yield the following result
\[\sum\limits_{n = 1}^\infty  {\frac{{{H_{n + \alpha }}}}{{\left( {n + k} \right)\left( {n + \alpha } \right)}}}  = \frac{{H_{\alpha  - k}^2 + H_{\alpha  - k}^{\left( 2 \right)}}}{{\alpha  - k}} - \sum\limits_{j = 1}^k {\frac{{{H_{\alpha  + j - k}}}}{{j\left( {\alpha  + j - k} \right)}}} \;\;\left( {\alpha  > k} \right).\tag{2.18}\]
Putting $f\left( {n,\alpha } \right) = {H_{n + \alpha }}$ in (1.10), and combining (2.18), we may easily deduce the desired result. \hfill$\square$
\begin{cor} For integers $r,k,m\in \N$ with $r\neq k$ and real $\alpha> max\{k,r\}$, we have
\begin{align*}
\sum\limits_{n = 1}^\infty  {\frac{{{H_{n + \alpha  - m}}}}{{\left( {n + r} \right)\left( {n + k} \right)}}}  =& \frac{{k - \alpha }}{{k - r}}\left\{ {\frac{{H_{\alpha  - k}^2 + H_{\alpha  - k}^{\left( 2 \right)}}}{{\alpha  - k}} - \sum\limits_{j = 1}^k {\frac{{{H_{\alpha  + j - k}}}}{{j\left( {\alpha  + j - k} \right)}}} } \right\}\\
&+ \frac{{\alpha  - r}}{{k - r}}\left\{ {\frac{{H_{\alpha  - r}^2 + H_{\alpha  - r}^{\left( 2 \right)}}}{{\alpha  - r}} - \sum\limits_{j = 1}^r {\frac{{{H_{\alpha  + j - r}}}}{{j\left( {\alpha  + j - r} \right)}}} } \right\}\\
&- \frac{1}{{k - r}}\sum\limits_{j = 1}^m {\left\{ {\frac{{{H_{j + \alpha  - m}} - {H_r}}}{{j + \alpha  - m - r}} - \frac{{{H_{j + \alpha  - m}} - {H_k}}}{{j + \alpha  - m - k}}} \right\}},\tag{2.19}
\end{align*}
where $\alpha  \ne m + r - j $and $m + k - j\;\left( {j = 1,2, \cdots ,m} \right)$ with $\alpha  - m \notin \N^-$.
\end{cor}
\pf By the definition of $H_{n+\alpha}$, we can find the following relation
\[{H_{n + \alpha }} - {H_{n + \alpha  - m}} = \sum\limits_{j = 1}^m {\frac{1}{{j + n + \alpha  - m}}} .\]
Hence, using the above identity, we can rewrite the series on the left hand side of (2.19) as
\[\sum\limits_{n = 1}^\infty  {\frac{{{H_{n + \alpha  - m}}}}{{\left( {n + r} \right)\left( {n + k} \right)}}}  = \sum\limits_{n = 1}^\infty  {\frac{{{H_{n + \alpha }}}}{{\left( {n + r} \right)\left( {n + k} \right)}}}  - \sum\limits_{j = 1}^m {\sum\limits_{n = 1}^\infty  {\frac{1}{{\left( {n + j + \alpha  - m} \right)\left( {n + r} \right)\left( {n + k} \right)}}} } .\tag{2.20}\]
By a direct calculation, we obtain
\begin{align*}
&\sum\limits_{n = 1}^\infty  {\frac{1}{{\left( {n + j + \alpha  - m} \right)\left( {n + r} \right)\left( {n + k} \right)}}} \\
&= \frac{1}{{k - r}}\left\{ {\sum\limits_{n = 1}^\infty  {\frac{1}{{\left( {n + j + \alpha  - m} \right)\left( {n + r} \right)}}}  - \sum\limits_{n = 1}^\infty  {\frac{1}{{\left( {n + j + \alpha  - m} \right)\left( {n + k} \right)}}} } \right\}\\
& = \frac{{{H_{j + \alpha  - m}} - {H_r}}}{{\left( {k - r} \right)\left( {j + \alpha  - m - r} \right)}} - \frac{{{H_{j + \alpha  - m}} - {H_k}}}{{\left( {k - r} \right)\left( {j + \alpha  - m - k} \right)}}.\tag{2.21}
\end{align*}
Substituting (2.16) and (2.21) into (2.20), we can prove (2.19). The proof of Corollary 2.3 is thus completed.\hfill$\square$\\
From (2.16) and (2.19), we have the following special cases:
\begin{align*}
&\sum\limits_{n = 1}^\infty  {\frac{{{H_{n - 1/2}}}}{{\left( {n + 1} \right)\left( {n + 2} \right)}}}  = \frac{2}{3} + \frac{1}{3}\ln 2,\ \sum\limits_{n = 1}^\infty  {\frac{{{H_{n + 1/2}}}}{{\left( {n + 1} \right)\left( {n + 2} \right)}}}  = 3\ln 2 - 1,\\
&\sum\limits_{n = 1}^\infty  {\frac{{{H_{n + 3/2}}}}{{\left( {n + 1} \right)\left( {n + 2} \right)}}}  = \frac{{14}}{3} - 5\ln 2,\ \sum\limits_{n = 1}^\infty  {\frac{{{H_{n + 5/2}}}}{{\left( {n + 1} \right)\left( {n + 2} \right)}}}  = \frac{{131}}{{45}} - \frac{7}{3}\ln 2.
\end{align*}
Next, we evaluate the summation of the shifted harmonic sums
of order two of the form,
\[\sum\limits_{n = 1}^\infty  {\frac{{H_{n + \alpha }^{\left( 2 \right)}}}{{\left( {n + r} \right)\left( {n + k} \right)}}} ,\;\;\alpha  > \max \{ k,r\} ,\;k,r \in \N_0,k \ne r.\]
First, we need to obtain the integral representation of ${H_{\alpha }^{\left( 2 \right)}}$. Setting $p=2$ in (1.11) gives
\[\int\limits_0^1 {{x^{\alpha  - 1}}{\rm{L}}{{\rm{i}}_2}\left( x \right)dx}  = \frac{{\zeta \left( 2 \right)}}{\alpha } - \frac{{{H_\alpha }}}{{{\alpha ^2}}}.\tag{2.22}\]
With the help of (2.13), we get
\[\frac{{H_\alpha ^{\left( 2 \right)}}}{\alpha } = \int\limits_0^1 {{x^{\alpha  - 1}}\ln x\ln \left( {1 - x} \right)} dx + \int\limits_0^1 {{x^{\alpha  - 1}}{\rm{L}}{{\rm{i}}_2}\left( x \right)dx},\ \alpha>0.\tag{2.23} \]
Replacing $\alpha$ by $n+\alpha$ in (2.23), then multiplying it by $(n+k)^{-1}$ and summing with respect to $n$, the result is
\begin{align*}
\sum\limits_{n = 1}^\infty  {\frac{{H_{n + \alpha }^{\left( 2 \right)}}}{{\left( {n + k} \right)\left( {n + \alpha } \right)}}}  =& \sum\limits_{n = 1}^\infty  {\frac{1}{{n + k}}\int\limits_0^1 {{x^{n + \alpha  - 1}}\ln x\ln \left( {1 - x} \right)} dx}  + \sum\limits_{n = 1}^\infty  {\frac{1}{{n + k}}\int\limits_0^1 {{x^{n + \alpha  - 1}}{\rm{L}}{{\rm{i}}_2}\left( x \right)dx} } \\
 = & - \int\limits_0^1 {{x^{\alpha  - k - 1}}\ln x{{\ln }^2}\left( {1 - x} \right)dx}  - \sum\limits_{j = 1}^k {\frac{1}{j}} \int\limits_0^1 {{x^{\alpha  + j - k - 1}}\ln x\ln \left( {1 - x} \right)} dx\\
 &- \int\limits_0^1 {{x^{\alpha  - k - 1}}\ln \left( {1 - x} \right){\rm{L}}{{\rm{i}}_2}\left( x \right)dx}  - \sum\limits_{j = 1}^k {\frac{1}{j}} \int\limits_0^1 {{x^{\alpha  + j - k - 1}}{\rm{L}}{{\rm{i}}_2}\left( x \right)} dx.\tag{2.24}
\end{align*}
We note that by using (2.11), the following identities ar easily derived
\begin{align*}
\int\limits_0^1 {{x^{\alpha  - 1}}\ln \left( {1 - x} \right){\rm{L}}{{\rm{i}}_2}\left( x \right)dx}  = & - \sum\limits_{n = 1}^\infty  {\frac{{{H_{n + \alpha }}}}{{{n^2}\left( {n + \alpha } \right)}}} \\
 =& \frac{1}{\alpha }\sum\limits_{n = 1}^\infty  {\frac{{{H_{n + \alpha }}}}{{n\left( {n + \alpha } \right)}}}  - \frac{1}{\alpha }\sum\limits_{n = 1}^\infty  {\frac{{{H_{n + \alpha }}}}{{{n^2}}}} .\tag{2.25}
\end{align*}
On the other hand, from (2.18), setting $k=0$, we have
\[\sum\limits_{n = 1}^\infty  {\frac{{{H_{n + \alpha }}}}{{n\left( {n + \alpha } \right)}}}  = \frac{{H_\alpha ^2 + H_\alpha ^{\left( 2 \right)}}}{\alpha },\alpha  \notin \N^-\cup\{0\}.\tag{2.26}\]
Hence, combining (2.13), (2.15), (2.22) and (2.24)-(2.26), we obtain
\begin{align*}
\sum\limits_{n = 1}^\infty  {\frac{{H_{n + \alpha }^{\left( 2 \right)}}}{{\left( {n + k} \right)\left( {n + \alpha } \right)}}}  =& \frac{1}{{\alpha  - k}}\sum\limits_{n = 1}^\infty  {\frac{{{H_{n + \alpha  - k}}}}{{{n^2}}}}  - 2\frac{{\zeta \left( 3 \right) - H_{\alpha  - k}^{\left( 3 \right)}}}{{\alpha  - k}}\\
 &- 2\frac{{\zeta \left( 2 \right) - H_{\alpha  - k}^{\left( 2 \right)}}}{{\alpha  - k}}{H_{\alpha  - k}} - \sum\limits_{j = 1}^k {\frac{{H_{\alpha  + j - k}^{\left( 2 \right)}}}{{j\left( {\alpha  + j - k} \right)}}} ,\ (\alpha>k) .\tag{2.27}
\end{align*}
Letting $f\left( {n,\alpha } \right) = H_{n + \alpha }^{\left( 2 \right)}$ in (1.10), then substituting (2.27) into (1.10), we get the following Theorem.
\begin{thm} For integers $r,k\geq 0$ and real $\alpha>k>r$, then we have
\begin{align*}
\sum\limits_{n = 1}^\infty  {\frac{{H_{n + \alpha }^{\left( 2 \right)}}}{{\left( {n + r} \right)\left( {n + k} \right)}}}  = \frac{1}{{k - r}}\left\{ \begin{array}{l}
 \left( {r - \alpha } \right)\sum\limits_{j = 1}^r {\frac{{H_{\alpha  + j - r}^{\left( 2 \right)}}}{{j\left( {\alpha  + j - r} \right)}}}  - \left( {k - \alpha } \right)\sum\limits_{j = 1}^k {\frac{{H_{\alpha  + j - k}^{\left( 2 \right)}}}{{j\left( {\alpha  + j - k} \right)}}}  \\
  - \sum\limits_{j = 1}^{k - r} {\frac{{{H_{\alpha  + j - k}}}}{{{{\left( {\alpha  + j - k} \right)}^2}}}}  + 2H_{\alpha  - r}^{\left( 3 \right)} + {H_{\alpha  - k}}\zeta \left( 2 \right) + 2{H_{\alpha  - r}}H_{\alpha  - r}^{\left( 2 \right)} \\
  - 2H_{\alpha  - k}^{\left( 3 \right)} - {H_{\alpha  - r}}\zeta \left( 2 \right) - 2{H_{\alpha  - k}}H_{\alpha  - k}^{\left( 2 \right)} \\
 \end{array} \right\}.\tag{2.28}
\end{align*}
\end{thm}
\pf From (1.10) and (2.27), we arrive at the conclusion that
\begin{align*}
\sum\limits_{n = 1}^\infty  \frac{{H_{n + \alpha }^{\left( 2 \right)}}}{{\left( {n + r} \right)\left( {n + k} \right)}} =&  \frac{1}{{k - r}}\sum\limits_{n = 1}^\infty  {\frac{{{H_{n + \alpha  - r}} - {H_{n + \alpha  - k}}}}{{{n^2}}}} \\
&+ \frac{{r - \alpha }}{{k - r}}\sum\limits_{j = 1}^r {\frac{{H_{\alpha  + j - r}^{\left( 2 \right)}}}{{j\left( {\alpha  + j - r} \right)}}}  - \frac{{k - \alpha }}{{k - r}}\sum\limits_{j = 1}^k {\frac{{H_{\alpha  + j - k}^{\left( 2 \right)}}}{{j\left( {\alpha  + j - k} \right)}}} \\
& + \frac{2}{{k - r}}\left\{ \begin{array}{l}
 H_{\alpha  - r}^{\left( 3 \right)} + {H_{\alpha  - k}}\zeta \left( 2 \right) + {H_{\alpha  - r}}H_{\alpha  - r}^{\left( 2 \right)} \\
  - H_{\alpha  - k}^{\left( 3 \right)} - {H_{\alpha  - r}}\zeta \left( 2 \right) - {H_{\alpha  - k}}H_{\alpha  - k}^{\left( 2 \right)} \\
 \end{array} \right\}.\tag{2.29}
\end{align*}
By the definition of shifted harmonic number, and using (1.14), we conclude that
\[H_{\alpha  - r}^{\left( m \right)} - H_{\alpha  - k}^{\left( m \right)} = \sum\limits_{j = 1}^{k - r} {\frac{1}{{{{\left( {j + \alpha  - k} \right)}^m}}}},\]
and for $2\leq p\in \N,\ 0\leq r<k\in \N$,
\begin{align*}
&\sum\limits_{n = 1}^\infty  {\frac{{{H_{n + \alpha  - r}} - {H_{n + \alpha  - k}}}}{{{n^p}}}}  =\sum\limits_{j = 1}^{k - r} {\sum\limits_{n = 1}^\infty  {\frac{1}{{{n^p}\left( {n + j + \alpha  - k} \right)}}} }\\
& = \sum\limits_{j = 1}^{k - r} {\sum\limits_{i = 1}^{p - 1} {\frac{{{{\left( { - 1} \right)}^{i - 1}}}}{{{{\left( {j + \alpha  - k} \right)}^i}}}\zeta \left( {p + 1 - i} \right)} }  + {\left( { - 1} \right)^{p - 1}}\sum\limits_{j = 1}^{k - r} { {\frac{{{H_{j + \alpha  - k}}}}{{{{\left( {j + \alpha  - k} \right)}^p}}}} } \\
 &= \sum\limits_{i = 1}^{p - 1} {{{\left( { - 1} \right)}^{i - 1}}\zeta \left( {p + 1 - i} \right)\left( {H_{\alpha  - r}^{\left( i \right)} - H_{\alpha  - k}^{\left( i \right)}} \right)}  + {\left( { - 1} \right)^{p - 1}}\sum\limits_{j = 1}^{k - r} { {\frac{{{H_{j + \alpha  - k}}}}{{{{\left( {j + \alpha  - k} \right)}^p}}}} } .\tag{2.30}
\end{align*}
Taking $p=2$ in (2.30) and substituting it into (2.29), by simple calculation, we obtain the desired result.\hfill$\square$
\begin{cor} For integers $r,k,m\in \N$ and $\alpha> k>r$, then
\begin{align*}
\sum\limits_{n = 1}^\infty  {\frac{{H_{n + \alpha  - m}^{\left( 2 \right)}}}{{\left( {n + r} \right)\left( {n + k} \right)}}}  =& \frac{1}{{k - r}}\left\{ \begin{array}{l}
 \left( {r - \alpha } \right)\sum\limits_{j = 1}^r {\frac{{H_{\alpha  + j - r}^{\left( 2 \right)}}}{{j\left( {\alpha  + j - r} \right)}}}  - \left( {k - \alpha } \right)\sum\limits_{j = 1}^k {\frac{{H_{\alpha  + j - k}^{\left( 2 \right)}}}{{j\left( {\alpha  + j - k} \right)}}}  \\
  - \sum\limits_{j = 1}^{k - r} {\frac{{{H_{\alpha  + j - k}}}}{{{{\left( {\alpha  + j - k} \right)}^2}}}}  + 2H_{\alpha  - r}^{\left( 3 \right)} + {H_{\alpha  - k}}\zeta \left( 2 \right) + 2{H_{\alpha  - r}}H_{\alpha  - r}^{\left( 2 \right)} \\
  - 2H_{\alpha  - k}^{\left( 3 \right)} - {H_{\alpha  - r}}\zeta \left( 2 \right) - 2{H_{\alpha  - k}}H_{\alpha  - k}^{\left( 2 \right)} \\
 \end{array} \right\}\\
 & - \frac{1}{{k - r}}\left\{ \begin{array}{l}
 \sum\limits_{j = 1}^m {\frac{{{H_{\alpha  + j - m}} - {H_r}}}{{{{\left( {\alpha  + j - m - r} \right)}^2}}}}  - \sum\limits_{j = 1}^m {\frac{{\zeta \left( 2 \right) - H_{\alpha  + j - m}^{\left( 2 \right)}}}{{\alpha  + j - m - r}}}  \\
  - \sum\limits_{j = 1}^m {\frac{{{H_{\alpha  + j - m}} - {H_k}}}{{{{\left( {\alpha  + j - m - k} \right)}^2}}}}  + \sum\limits_{j = 1}^m {\frac{{\zeta \left( 2 \right) - H_{\alpha  + j - m}^{\left( 2 \right)}}}{{\alpha  + j - m - k}}}  \\
 \end{array} \right\}.\tag{2.31}
\end{align*}
where $\alpha  \ne m + r - j,m + k - j\;\left( {j = 1,2, \cdots ,m} \right)$ and $\alpha  - m \notin \N^-$.
\end{cor}
\pf By a similar argument as in the proof of Corollary 2.3, we have
\[\sum\limits_{n = 1}^\infty  {\frac{{H_{n + \alpha  - m}^{\left( 2 \right)}}}{{\left( {n + r} \right)\left( {n + k} \right)}}}  = \sum\limits_{n = 1}^\infty  {\frac{{H_{n + \alpha }^{\left( 2 \right)}}}{{\left( {n + r} \right)\left( {n + k} \right)}}}  - \sum\limits_{j = 1}^m {\sum\limits_{n = 1}^\infty  {\frac{1}{{{{\left( {n + j + \alpha  - m} \right)}^2}\left( {n + k} \right)\left( {n + r} \right)}}} }.\tag{2.32} \]
On the other hand, we easily obtain the result
\begin{align*}
&\sum\limits_{n = 1}^\infty  {\frac{1}{{{{\left( {n + j + \alpha  - m} \right)}^2}\left( {n + k} \right)\left( {n + r} \right)}}} \\
& = \frac{1}{{k - r}}\sum\limits_{n = 1}^\infty  {\left( {\frac{1}{{{{\left( {n + j + \alpha  - m} \right)}^2}\left( {n + r} \right)}} - \frac{1}{{{{\left( {n + j + \alpha  - m} \right)}^2}\left( {n + k} \right)}}} \right)} \\
& = \frac{1}{{k - r}}\left\{ \begin{array}{l}
 \frac{1}{{{{\left( {j + \alpha  - m - r} \right)}^2}}}\sum\limits_{n = 1}^\infty  {\left( {\frac{1}{{n + r}} - \frac{1}{{n + j + \alpha  - m}}} \right)}  \\
  - \frac{1}{{j + \alpha  - m - r}}\sum\limits_{n = 1}^\infty  {\frac{1}{{{{\left( {n + j + \alpha  - m} \right)}^2}}}}  \\
  - \frac{1}{{{{\left( {j + \alpha  - m - k} \right)}^2}}}\sum\limits_{n = 1}^\infty  {\left( {\frac{1}{{n + k}} - \frac{1}{{n + j + \alpha  - m}}} \right)}  \\
  + \frac{1}{{j + \alpha  - m - k}}\sum\limits_{n = 1}^\infty  {\frac{1}{{{{\left( {n + j + \alpha  - m} \right)}^2}}}}  \\
 \end{array} \right\}\\
 & = \frac{1}{{k - r}}\left\{ \begin{array}{l}
 \sum\limits_{j = 1}^m {\frac{{{H_{\alpha  + j - m}} - {H_r}}}{{{{\left( {\alpha  + j - m - r} \right)}^2}}}}  - \sum\limits_{j = 1}^m {\frac{{\zeta \left( 2 \right) - H_{\alpha  + j - m}^{\left( 2 \right)}}}{{\alpha  + j - m - r}}}  \\
  - \sum\limits_{j = 1}^m {\frac{{{H_{\alpha  + j - m}} - {H_k}}}{{{{\left( {\alpha  + j - m - k} \right)}^2}}}}  + \sum\limits_{j = 1}^m {\frac{{\zeta \left( 2 \right) - H_{\alpha  + j - m}^{\left( 2 \right)}}}{{\alpha  + j - m - k}}}  \\
 \end{array} \right\}.\tag{2.33}
\end{align*}
Substituting (2.28) and (2.33) into (2.32) respectively, we deduce (2.31). This completes the proof of Corollary 2.5.\hfill$\square$
\begin{thm} For integers $r,k\geq 0$ and $\alpha>k>r$, then we have
\begin{align*}
\sum\limits_{n = 1}^\infty  {\frac{{H_{n + \alpha }^2 + H_{n + \alpha }^{\left( 2 \right)}}}{{\left( {n + r} \right)\left( {n + k} \right)}}}  =& \frac{{k - \alpha }}{{k - r}}\left\{ {\frac{{H_{\alpha  - k}^3 + 3{H_{\alpha  - k}}H_{\alpha  - k}^{\left( 2 \right)} + 2H_{\alpha  - k}^{\left( 3 \right)}}}{{\alpha  - k}} - \sum\limits_{j = 1}^k {\frac{{H_{\alpha  + j - k}^2 + H_{\alpha  + j - k}^{\left( 2 \right)}}}{{j\left( {\alpha  + j - k} \right)}}} } \right\}\\
& + \frac{{\alpha  - r}}{{k - r}}\left\{ {\frac{{H_{\alpha  - r}^3 + 3{H_{\alpha  - r}}H_{\alpha  - r}^{\left( 2 \right)} + 2H_{\alpha  - r}^{\left( 3 \right)}}}{{\alpha  - r}} - \sum\limits_{j = 1}^r {\frac{{H_{\alpha  + j - r}^2 + H_{\alpha  + j - r}^{\left( 2 \right)}}}{{j\left( {\alpha  + j - r} \right)}}} } \right\}.\tag{2.34}
\end{align*}
\end{thm}
\pf From (2.12) and (2.14), we know that
\begin{align*}
\sum\limits_{n = 1}^\infty  {\frac{{H_{n + \alpha }^2 + H_{n + \alpha }^{\left( 2 \right)}}}{{\left( {n + k} \right)\left( {n + \alpha } \right)}}}  &=\sum\limits_{n = 1}^\infty  {\frac{1}{{n + k}}\int\limits_0^1 {{x^{n + \alpha  - 1}}{{\ln }^2}\left( {1 - x} \right)} } dx\\
& =  - \int\limits_0^1 {{x^{\alpha  - k - 1}}{{\ln }^3}\left( {1 - x} \right)} dx - \sum\limits_{j = 1}^k {\frac{1}{j}\int\limits_0^1 {{x^{\alpha  + j - k - 1}}{{\ln }^2}\left( {1 - x} \right)} dx} \\
& = \frac{{H_{\alpha  - k}^3 + 3{H_{\alpha  - k}}H_{\alpha  - k}^{\left( 2 \right)} + 2H_{\alpha  - k}^{\left( 3 \right)}}}{{\alpha  - k}} - \sum\limits_{j = 1}^k {\frac{{H_{\alpha  + j - k}^2 + H_{\alpha  + j - k}^{\left( 2 \right)}}}{{j\left( {\alpha  + j - k} \right)}}} .\tag{2.35}
\end{align*}
Setting $f\left( {n,\alpha } \right) = H_{n + \alpha }^2 + H_{n + \alpha }^{\left( 2 \right)}$ in (1.10) and combining (2.35), the result is (2.34).\hfill$\square$\\
Similarly to the proof of Theorem 2.6, we have the following similar result.
\begin{thm} For integers $r,k\geq 0$ and $\alpha>k>r$, then the following identity holds:
\begin{align*}
&\sum\limits_{n = 1}^\infty  {\frac{{H_{n + \alpha }^3 + 3{H_{n + \alpha }}H_{n + \alpha }^{\left( 2 \right)} + 2H_{n + \alpha }^{\left( 3 \right)}}}
{{\left( {n + r} \right)\left( {n + k} \right)}}} \\
&  = \frac{{k - \alpha }}
{{k - r}}\left\{ \begin{gathered}
  \frac{{H_{\alpha  - k}^4 + 6H_{\alpha  - k}^2H_{\alpha  - k}^{\left( 2 \right)} + 8{H_{\alpha  - k}}H_{\alpha  - k}^{\left( 3 \right)} + 3{{\left( {H_{\alpha  - k}^{\left( 2 \right)}} \right)}^2} + 6H_{\alpha  - k}^{\left( 4 \right)}}}
{{\alpha  - k}} \hfill \\
   - \sum\limits_{j = 1}^k {\frac{{H_{\alpha  + j - k}^3 + 3{H_{\alpha  + j - k}}H_{\alpha  + j - k}^{\left( 2 \right)} + 2H_{\alpha  + j - k}^{\left( 3 \right)}}}
{{j\left( {\alpha  + j - k} \right)}}}  \hfill \\
\end{gathered}  \right\}\\
& \quad + \frac{{\alpha  - r}}
{{k - r}}\left\{ \begin{gathered}
  \frac{{H_{\alpha  - r}^4 + 6H_{\alpha  - r}^2H_{\alpha  - r}^{\left( 2 \right)} + 8{H_{\alpha  - r}}H_{\alpha  - r}^{\left( 3 \right)} + 3{{\left( {H_{\alpha  - r}^{\left( 2 \right)}} \right)}^2} + 6H_{\alpha  - r}^{\left( 4 \right)}}}
{{\alpha  - r}} \hfill \\
   - \sum\limits_{j = 1}^r {\frac{{H_{\alpha  + j - r}^3 + 3{H_{\alpha  + j - r}}H_{\alpha  + j - r}^{\left( 2 \right)} + 2H_{\alpha  + j - r}^{\left( 3 \right)}}}
{{j\left( {\alpha  + j - r} \right)}}}  \hfill \\
\end{gathered}  \right\} .\tag{2.36}
\end{align*}
\end{thm}
\pf Applying the same arguments as in the proof of Theorem 2.6, we may easily deduce
\begin{align*}
&\sum\limits_{n = 1}^\infty  {\frac{{H_{n + \alpha }^3 + 3{H_{n + \alpha }}H_{n + \alpha }^{\left( 2 \right)} + 2H_{n + \alpha }^{\left( 3 \right)}}}
{{\left( {n + \alpha } \right)\left( {n + k} \right)}}} \\ &= \int\limits_0^1 {{x^{\alpha  - k - 1}}{{\ln }^4}\left( {1 - x} \right)} dx + \sum\limits_{j = 1}^k {\frac{1}
{j}\int\limits_0^1 {{x^{\alpha  + j - k - 1}}{{\ln }^3}\left( {1 - x} \right)} dx} .\tag{2.37}
\end{align*}
Setting $m=4$ in (2.7) and combining (1.10), (2.11), (2.12) with (2.14) we obtain
 \[\int\limits_0^1 {{x^{\alpha  - 1}}{{\ln }^4}\left( {1 - x} \right)} dx = \frac{{H_\alpha ^4 + 6H_\alpha ^2H_\alpha ^{\left( 2 \right)} + 8{H_\alpha }H_\alpha ^{\left( 3 \right)} + 3{{\left( {H_\alpha ^{\left( 2 \right)}} \right)}^2} + 6H_\alpha ^{\left( 4 \right)}}}
{\alpha }.\tag{2.38}\]
Therefore, by using (2.14), (2.37) and (2.38) yields the desired result. This completes the proof Theorem 2.7. \hfill$\square$\\
Substituting (2.28) into (2.34), we can obtain the following Corollary.
\begin{cor} For integers $r,k\geq 0$ and $\alpha>k>r$, then we have the quadratic sums
\begin{align*}
\sum\limits_{n = 1}^\infty  {\frac{{H_{n + \alpha }^2}}{{\left( {n + r} \right)\left( {n + k} \right)}}}  = \frac{1}{{k - r}}\left\{ \begin{array}{l}
 H_{\alpha  - r}^3 + {H_{\alpha  - r}}H_{\alpha  - r}^{\left( 2 \right)} + {H_{\alpha  - r}}\zeta \left( 2 \right) \\
  - H_{\alpha  - k}^3 - {H_{\alpha  - k}}H_{\alpha  - k}^{\left( 2 \right)} - {H_{\alpha  - k}}\zeta \left( 2 \right) \\
  + \left( {r - \alpha } \right)\sum\limits_{j = 1}^r {\frac{{H_{\alpha  + j - r}^2}}{{j\left( {\alpha  + j - r} \right)}}}  + \sum\limits_{j = 1}^{k - r} {\frac{{{H_{\alpha  + j - k}}}}{{{{\left( {\alpha  + j - k} \right)}^2}}}}  \\
  - \left( {k - \alpha } \right)\sum\limits_{j = 1}^k {\frac{{H_{\alpha  + j - k}^2}}{{j\left( {\alpha  + j - k} \right)}}}  \\
 \end{array} \right\}.\tag{2.39}
 \end{align*}
\end{cor}
Further, using the definition of shifted harmonic number, by a simple calculation, the following relations are easily derived
\begin{align*}
&\frac{\partial }{{\partial \alpha }}\left( {H_{n + \alpha }^{\left( m \right)}} \right) = m\left( {\zeta \left( {m + 1} \right) - H_{n + \alpha }^{\left( {m + 1} \right)}} \right),\\
&\frac{{{\partial ^m}}}{{\partial {\alpha ^m}}}\left( {{H_{n + \alpha }}} \right) = {\left( { - 1} \right)^{m + 1}}m!\left( {\zeta \left( {m + 1} \right) - H_{n + \alpha }^{\left( {m + 1} \right)}} \right),\\
&\frac{\partial }{{\partial \alpha }}\left( {H_{n + \alpha }^2} \right) = 2\zeta \left( 2 \right){H_{n + \alpha }} - 2{H_{n + \alpha }}H_{n + \alpha }^{\left( 2 \right)},\\
&\frac{\partial }{{\partial \alpha }}\left( {H_{n + \alpha }^3} \right) = 3\zeta \left( 2 \right)H_{n + \alpha }^2 - 3H_{n + \alpha }^2H_{n + \alpha }^{\left( 2 \right)}.
\end{align*}
Hence, from Theorem 2.2, 2.4, 2.7 and Corollary 2.8 with the above relations, we can get the Theorem 2.9.
\begin{thm} For positive integers $k,r,m$ and real $\alpha>{\rm max}\{k,r\}$ with $k\neq r$. Then the linear, quadratic and cubic sums
\[\sum\limits_{n = 1}^\infty  {\frac{{H_{n + \alpha }^{\left( m \right)}}}{{\left( {n + r} \right)\left( {n + k} \right)}}} ,\sum\limits_{n = 1}^\infty  {\frac{{{H_{n + \alpha }}H_{n + \alpha }^{\left( 2 \right)}}}{{\left( {n + r} \right)\left( {n + k} \right)}}} ,\sum\limits_{n = 1}^\infty  {\frac{{H_{n + \alpha }^3}}{{\left( {n + r} \right)\left( {n + k} \right)}}} ,\sum\limits_{n = 1}^\infty  {\frac{{H_{n + \alpha }^2H_{n + \alpha }^{\left( 2 \right)}}}{{\left( {n + r} \right)\left( {n + k} \right)}}} \]
can be expressed in terms of shifted harmonic numbers and zeta values.
\end{thm}
As simple example is as follows:
\begin{cor} For integers $r,k\geq 0$ and $\alpha>k>r$, then we have
\[\sum\limits_{n = 1}^\infty  {\frac{{H_{n + \alpha }^{\left( 3 \right)}}}{{\left( {n + r} \right)\left( {n + k} \right)}}}  = \frac{1}{{k - r}}\left\{ \begin{array}{l}
 \left( {r - \alpha } \right)\sum\limits_{j = 1}^r {\frac{{H_{\alpha  + j - r}^{\left( 3 \right)}}}{{j\left( {\alpha  + j - r} \right)}}}  - \left( {k - \alpha } \right)\sum\limits_{j = 1}^k {\frac{{H_{\alpha  + j - k}^{\left( 3 \right)}}}{{j\left( {\alpha  + j - k} \right)}}}  \\
  - \frac{1}{2}\sum\limits_{j = 1}^{k - r} {\frac{{H_{\alpha  + j - k}^{\left( 2 \right)}}}{{{{\left( {\alpha  + j - k} \right)}^2}}}}  - \sum\limits_{j = 1}^{k - r} {\frac{{{H_{\alpha  + j - k}}}}{{{{\left( {\alpha  + j - k} \right)}^3}}}}  + 3\left( {H_{\alpha  - r}^{\left( 4 \right)} - H_{\alpha  - k}^{\left( 4 \right)}} \right) \\
  + \left( {H_{\alpha  - k}^{\left( 2 \right)} - H_{\alpha  - r}^{\left( 2 \right)}} \right)\zeta \left( 2 \right) + \left( {{H_{\alpha  - k}} - {H_{\alpha  - r}}} \right)\zeta \left( 3 \right) \\
  + \left( {{{\left( {H_{\alpha  - r}^{\left( 2 \right)}} \right)}^2} - {{\left( {H_{\alpha  - k}^{\left( 2 \right)}} \right)}^2}} \right) + 2\left( {{H_{\alpha  - r}}H_{\alpha  - r}^{\left( 3 \right)} - {H_{\alpha  - k}}H_{\alpha  - k}^{\left( 3 \right)}} \right) \\
 \end{array} \right\}.\]
\end{cor}
\section{Some results of ${W^{(l)}_{k,r}}\left( {p,m,\alpha } \right)$ }
In this section, we give some closed form form sums of ${W^{(l)}_{k,r}}\left( {p,m,\alpha } \right)$ through shifted harmonic numbers and zeta values. First, we consider the partial fraction decomposition
\[\frac{1}{{\prod\limits_{i = 1}^m {\left( {n + {a_i}} \right)} }} = \sum\limits_{j = 1}^m {\frac{{{A_j}}}{{n + {a_j}}}},\tag{3.1} \]
where
\[{A_j} = \mathop {\lim }\limits_{n \to  - {a_j}} \frac{{n + {a_j}}}{{\prod\limits_{i = 1}^m {\left( {n + {a_i}} \right)} }} = \prod\limits_{i = 1,i \ne j}^n {{{\left( {{a_i} - {a_j}} \right)}^{ - 1}}}.\]
Letting $m=k$ and $a_i=r+i$ in (3.1), we have
\[\frac{1}{{\left( {\begin{array}{*{20}{c}}
   {n + k + r}  \\
   k  \\
\end{array}} \right)}} = \frac{{k!}}{{\prod\limits_{i = 1}^k {\left( {n + r + i} \right)} }} = \sum\limits_{j = 1}^k {{{\left( { - 1} \right)}^{j + 1}}j\left( {\begin{array}{*{20}{c}}
   k  \\
   j  \\
\end{array}} \right)} \frac{1}{{n + r + j}},\ k,r\in \N_0.\tag{3.2}\]
On the other hand, we can find that
\begin{align*}
\frac{1}{{\left( {\begin{array}{*{20}{c}}
   {n + k + r}  \\
   k  \\
\end{array}} \right)}} =& \frac{k}{{\left( {n + r + 1} \right)\left( {\begin{array}{*{20}{c}}
   {n + k + r}  \\
   {k - 1}  \\
\end{array}} \right)}}\\
 = &\frac{k}{{\left( {n + r + 1} \right)}}\sum\limits_{j = 1}^{k - 1} {{{\left( { - 1} \right)}^{j + 1}}j\left( {\begin{array}{*{20}{c}}
   {k - 1}  \\
   j  \\
\end{array}} \right)} \frac{1}{{n + r + 1 + j}}\\
 =& k\sum\limits_{j = 1}^{k - 1} {{{\left( { - 1} \right)}^{j + 1}}j\left( {\begin{array}{*{20}{c}}
   {k - 1}  \\
   j  \\
\end{array}} \right)} \frac{1}{{\left( {n + r + 1} \right)\left( {n + r + 1 + j} \right)}},\ r\in\N_0,k\in\N.\tag{3.3}
\end{align*}
From identities (2.16), (2.28), (2.39), (3.2) and (3.3), we obtain the following new results: for $\alpha>k+r,\ 2\leq k\in \N,\ r\in \N_0$,
\begin{align*}
W_{k,r}^{\left( 1 \right)}\left( {0,1,\alpha } \right) &= \sum\limits_{n = 1}^\infty  {\frac{{{H_{n + \alpha }}}}{{\left( {\begin{array}{*{20}{c}}
   {n + k + r}  \\
   k  \\
\end{array}} \right)}}} \\
& = k\sum\limits_{j = 1}^{k - 1} {{{\left( { - 1} \right)}^{j + 1}}j\left( {\begin{array}{*{20}{c}}
   {k - 1}  \\
   j  \\
\end{array}} \right)} \sum\limits_{n = 1}^\infty  {\frac{{{H_{n + \alpha }}}}{{\left( {n + r + 1} \right)\left( {n + r + 1 + j} \right)}}} \\
& = k\sum\limits_{j = 1}^{k - 1} {{{\left( { - 1} \right)}^{j + 1}}\left( {\begin{array}{*{20}{c}}
   {k - 1}  \\
   j  \\
\end{array}} \right)\left\{ \begin{array}{l}
 H_{\alpha  - r - 1}^2 + H_{\alpha  - r - 1}^{\left( 2 \right)} - H_{\alpha  - r - 1 - j}^2 - H_{\alpha  - r - 1 - j}^{\left( 2 \right)} \\
  - \left( {r + 1 + j - \alpha } \right)\sum\limits_{i = 1}^{r + 1 + j} {\frac{{{H_{\alpha  + i - r - 1 - j}}}}{{i\left( {\alpha  + i - r - 1 - j} \right)}}}  \\
  + \left( {r + 1 - \alpha } \right)\sum\limits_{i = 1}^{r + 1} {\frac{{{H_{\alpha  + i - r - 1}}}}{{i\left( {\alpha  + i - r - 1} \right)}}}  \\
 \end{array} \right\}}.\tag{3.4}\\
W_{k,r}^{\left( 1 \right)}\left( {0,2,\alpha } \right) &= \sum\limits_{n = 1}^\infty  {\frac{{H_{n + \alpha }^{\left( 2 \right)}}}{{\left( {\begin{array}{*{20}{c}}
   {n + k + r}  \\
   k  \\
\end{array}} \right)}}} \\
& = k\sum\limits_{j = 1}^{k - 1} {{{\left( { - 1} \right)}^{j + 1}}j\left( {\begin{array}{*{20}{c}}
   {k - 1}  \\
   j  \\
\end{array}} \right)} \sum\limits_{n = 1}^\infty  {\frac{{H_{n + \alpha }^{\left( 2 \right)}}}{{\left( {n + r + 1} \right)\left( {n + r + 1 + j} \right)}}} \\
& = k\sum\limits_{j = 1}^{k - 1} {{{\left( { - 1} \right)}^{j + 1}}\left( {\begin{array}{*{20}{c}}
   {k - 1}  \\
   j  \\
\end{array}} \right)\left\{ \begin{array}{l}
 \left( {r + 1 - \alpha } \right)\sum\limits_{i = 1}^{r + 1} {\frac{{H_{\alpha  + i - r - 1}^{\left( 2 \right)}}}{{i\left( {\alpha  + i - r - 1} \right)}}}  \\
  - \left( {r + 1 + j - \alpha } \right)\sum\limits_{i = 1}^{r+1+j} {\frac{{H_{\alpha  + i - r - 1 - j}^{\left( 2 \right)}}}{{i\left( {\alpha  + i - r - 1 - j} \right)}}}  \\
  - \sum\limits_{i = 1}^{j} {\frac{{{H_{\alpha  + i - r - 1 - j}}}}{{{{\left( {\alpha  + i - r - 1 - j} \right)}^2}}}}  \\
  + 2H_{\alpha  - r - 1}^{\left( 3 \right)} + {H_{\alpha  - r - 1 - j}}\zeta \left( 2 \right) + 2{H_{\alpha  - r - 1}}H_{\alpha  - r - 1}^{\left( 2 \right)} \\
  - 2H_{\alpha  - r - 1 - j}^{\left( 3 \right)} - {H_{\alpha  - r - 1}}\zeta \left( 2 \right) - 2{H_{\alpha  - r - 1 - j}}H_{\alpha  - r - 1 - j}^{\left( 2 \right)} \\
 \end{array} \right\}},\tag{3.5}\\
W_{k,r}^{\left( 2 \right)}\left( {0,1,\alpha } \right) &= \sum\limits_{n = 1}^\infty  {\frac{{H_{n + \alpha }^2}}{{\left( {\begin{array}{*{20}{c}}
   {n + k + r}  \\
   k  \\
\end{array}} \right)}}} \\
& = k\sum\limits_{j = 1}^{k - 1} {{{\left( { - 1} \right)}^{j + 1}}j\left( {\begin{array}{*{20}{c}}
   {k - 1}  \\
   j  \\
\end{array}} \right)} \sum\limits_{n = 1}^\infty  {\frac{{H_{n + \alpha }^2}}{{\left( {n + r + 1} \right)\left( {n + r + 1 + j} \right)}}} \\
& = k\sum\limits_{j = 1}^{k - 1} {{{\left( { - 1} \right)}^{j + 1}}\left( {\begin{array}{*{20}{c}}
   {k - 1}  \\
   j  \\
\end{array}} \right)} \left\{ \begin{array}{l}
 H_{\alpha  - r - 1}^3 + {H_{\alpha  - r - 1}}H_{\alpha  - r - 1}^{\left( 2 \right)} + {H_{\alpha  - r - 1}}\zeta \left( 2 \right) \\
  - H_{\alpha  - r - 1 - j}^3 - {H_{\alpha  - r - 1 - j}}H_{\alpha  - r - 1 - j}^{\left( 2 \right)} - {H_{\alpha  - r - 1 - j}}\zeta \left( 2 \right) \\
  + \left( {r + 1 - \alpha } \right)\sum\limits_{i = 1}^{r + 1} {\frac{{H_{\alpha  + i - r - 1}^2}}{{i\left( {\alpha  + i - r - 1} \right)}}}  \\
  + \sum\limits_{i = 1}^j {\frac{{{H_{\alpha  + i - r - 1 - j}}}}{{{{\left( {\alpha  + i - r - 1 - j} \right)}^2}}}}  \\
  - \left( {r + 1 + j - \alpha } \right)\sum\limits_{i = 1}^{r + 1 + j} {\frac{{H_{\alpha  + i - r - 1 - j}^2}}{{i\left( {\alpha  + i - r - 1 - j} \right)}}}  \\
 \end{array} \right\}.\tag{3.6}
\end{align*}
and for $k\in \N$ and $\alpha>k$,
\begin{align*}
W_{k,0}^{\left( 1 \right)}\left( {1,1,\alpha } \right) &= \sum\limits_{n = 1}^\infty  {\frac{{{H_{n + \alpha }}}}{{n\left( {\begin{array}{*{20}{c}}
   {n + k}  \\
   k  \\
\end{array}} \right)}}} \\
& = \sum\limits_{j = 1}^k {{{\left( { - 1} \right)}^{j + 1}}j\left( {\begin{array}{*{20}{c}}
   k  \\
   j  \\
\end{array}} \right)} \sum\limits_{n = 1}^\infty  {\frac{{{H_{n + \alpha }}}}{{n\left( {n + j} \right)}}} \\
& = \sum\limits_{j = 1}^k {{{\left( { - 1} \right)}^{j + 1}}\left( {\begin{array}{*{20}{c}}
   k  \\
   j  \\
\end{array}} \right)} \left\{ \begin{array}{l}
 H_\alpha ^2 + H_\alpha ^{\left( 2 \right)} - H_{\alpha  - j}^2 - H_{\alpha  - j}^{\left( 2 \right)} \\
  - \left( {j - \alpha } \right)\sum\limits_{i = 1}^j {\frac{{{H_{\alpha  + i - j}}}}{{i\left( {\alpha  + i - j} \right)}}}  \\
 \end{array} \right\},\tag{3.7}\\
W_{k,0}^{\left( 1 \right)}\left( {1,2,\alpha } \right) &= \sum\limits_{n = 1}^\infty  {\frac{{H_{n + \alpha }^{\left( 2 \right)}}}{{n\left( {\begin{array}{*{20}{c}}
   {n + k}  \\
   k  \\
\end{array}} \right)}}} \\
& = \sum\limits_{j = 1}^k {{{\left( { - 1} \right)}^{j + 1}}j\left( {\begin{array}{*{20}{c}}
   k  \\
   j  \\
\end{array}} \right)} \sum\limits_{n = 1}^\infty  {\frac{{{H^{(2)}_{n + \alpha }}}}{{n\left( {n + j} \right)}}} \\
& = \sum\limits_{j = 1}^k {{{\left( { - 1} \right)}^{j + 1}}\left( {\begin{array}{*{20}{c}}
   k  \\
   j  \\
\end{array}} \right)\left\{ \begin{array}{l}
 2H_\alpha ^{\left( 3 \right)} + {H_{\alpha  - j}}\zeta \left( 2 \right) + 2{H_\alpha }H_\alpha ^{\left( 2 \right)} \\
  - 2H_{\alpha  - j}^{\left( 3 \right)} - {H_\alpha }\zeta \left( 2 \right) - 2{H_{\alpha  - j}}H_{\alpha  - j}^{\left( 2 \right)} \\
  - \left( {j - \alpha } \right)\sum\limits_{i = 1}^j {\frac{{H_{\alpha  + i - j}^{\left( 2 \right)}}}{{i\left( {\alpha  + i - j} \right)}}}  \\
  - \sum\limits_{i = 1}^j {\frac{{{H_{\alpha  + i - j}}}}{{{{\left( {\alpha  + i - j} \right)}^2}}}}  \\
 \end{array} \right\}},\tag{3.8} \\
W_{k,0}^{\left( 2 \right)}\left( {1,1,\alpha } \right) &= \sum\limits_{n = 1}^\infty  {\frac{{H_{n + \alpha }^2}}{{n\left( {\begin{array}{*{20}{c}}
   {n + k}  \\
   k  \\
\end{array}} \right)}}} \\
& = \sum\limits_{j = 1}^k {{{\left( { - 1} \right)}^{j + 1}}j\left( {\begin{array}{*{20}{c}}
   k  \\
   j  \\
\end{array}} \right)} \sum\limits_{n = 1}^\infty  {\frac{{H_{n + \alpha }^2}}{{n\left( {n + j} \right)}}} \\
& = \sum\limits_{j = 1}^k {{{\left( { - 1} \right)}^{j + 1}}\left( {\begin{array}{*{20}{c}}
   k  \\
   j  \\
\end{array}} \right)\left\{ \begin{array}{l}
 H_\alpha ^3 + {H_\alpha }H_\alpha ^{\left( 2 \right)} + {H_\alpha }\zeta \left( 2 \right) \\
  - H_{\alpha  - j}^3 - {H_{\alpha  - j}}H_{\alpha  - j}^{\left( 2 \right)} - {H_{\alpha  - j}}\zeta \left( 2 \right) \\
  + \sum\limits_{i = 1}^j {\frac{{{H_{\alpha  + i - j}}}}{{{{\left( {\alpha  + i - j} \right)}^2}}}}  \\
  - \left( {j - \alpha } \right)\sum\limits_{i = 1}^j {\frac{{H_{\alpha  + i - j}^2}}{{i\left( {\alpha  + i - j} \right)}}}  \\
 \end{array} \right\}}.\tag{3.9}
\end{align*}
Hence, from Corollary 2.8, Theorem 2.9 and formulas (3.2), (3.3), we obtain the following description of ${W^{(l)}_{k,r}}\left( {p,m,\alpha } \right)$.
\begin{thm} For positive integers $k,r,m$ and real $\alpha\ (\alpha>k>r)$ with $p=0,1$, then the sums
\[W_{k,r}^{(1)}\left( {p,m,\alpha } \right),W_{k,r}^{(2)}\left( {p,1,\alpha } \right),W_{k,r}^{(3)}\left( {p,1,\alpha } \right)\]
and \[\sum\limits_{n = 1}^\infty  {\frac{{{H_{n + \alpha }}H_{n + \alpha }^{\left( 2 \right)}}}{{{n^p}\left( {\begin{array}{*{20}{c}}
   {n + k + r}  \\
   k  \\
\end{array}} \right)}}} ,\sum\limits_{n = 1}^\infty  {\frac{{H_{n + \alpha }^2H_{n + \alpha }^{\left( 2 \right)}}}{{{n^p}\left( {\begin{array}{*{20}{c}}
   {n + k + r}  \\
   k  \\
\end{array}} \right)}}} \]
can be expressed in terms of shifted harmonic numbers and ordinary zeta values.
\end{thm}
At the end of this section we give a explicit formula of sums associated with shifted harmonic numbers. We define the parametric polylogarithm function by the series
\[{\rm{L}}{{\rm{i}}_{p,\alpha }}\left( x \right): = \sum\limits_{n = 1}^\infty  {\frac{{{x^n}}}{{{{\left( {n + \alpha } \right)}^p}}}} ,\;x \in \left( { - 1,1} \right),\;\Re(p)>1,\ \alpha \notin  \N^-.\]
Next, we consider the following integral
\[\int\limits_0^1 {{x^{r - 1}}{\rm{L}}{{\rm{i}}_{p,\alpha }}\left( x \right){\rm{L}}{{\rm{i}}_{m,\beta }}\left( x \right)} dx,\;p,m\in \N,\ \alpha,\beta,r \notin \N^-.\]
First, using integration by parts, the following identity is easily derived
\[\int\limits_0^1 {{x^{r - 1}}{\rm{L}}{{\rm{i}}_{p,\alpha }}\left( x \right)} dx = \sum\limits_{i = 1}^{p - 1} {\frac{{{{\left( { - 1} \right)}^{i - 1}}}}{{{{\left( {r - \alpha } \right)}^i}}}\zeta \left( {p + 1 - i,\alpha  + 1} \right)}  + {\left( { - 1} \right)^{p - 1}}\frac{{{H_r} - {H_\alpha }}}{{{{\left( {r - \alpha } \right)}^p}}}.\tag{3.10}\]
We note that
\begin{align*}
\int\limits_0^1 {{x^{r - 1}}{\rm{L}}{{\rm{i}}_{p,\alpha }}\left( x \right){\rm{L}}{{\rm{i}}_{m,\beta }}\left( x \right)} dx &= \sum\limits_{n = 1}^\infty  {\frac{1}{{{{\left( {n + \alpha } \right)}^p}}}\int\limits_0^1 {{x^{n+r - 1}}{\rm{L}}{{\rm{i}}_{m,\beta }}\left( x \right)} dx} \\
& = \sum\limits_{n = 1}^\infty  {\frac{1}{{{{\left( {n + \beta } \right)}^m}}}\int\limits_0^1 {{x^{n+r - 1}}{\rm{L}}{{\rm{i}}_{p,\alpha }}\left( x \right)} dx} .\tag{3.11}
\end{align*}
Substituting (3.11) into (3.10) respectively, we can deduce that
\begin{align*}
&{\left( { - 1} \right)^{m - 1}}\sum\limits_{n = 1}^\infty  {\frac{{{H_{n + r}}}}{{{{\left( {n + \alpha } \right)}^p}{{\left( {n + r - \beta } \right)}^m}}}}  - {\left( { - 1} \right)^{p - 1}}\sum\limits_{n = 1}^\infty  {\frac{{{H_{n + r}}}}{{{{\left( {n + \beta } \right)}^m}{{\left( {n + r - \alpha } \right)}^p}}}} \\
& = {\left( { - 1} \right)^{m - 1}}{H_\beta }\sum\limits_{n = 1}^\infty  {\frac{1}{{{{\left( {n + \alpha } \right)}^p}{{\left( {n + r - \beta } \right)}^m}}}}  - {\left( { - 1} \right)^{p - 1}}{H_\alpha }\sum\limits_{n = 1}^\infty  {\frac{1}{{{{\left( {n + \beta } \right)}^m}{{\left( {n + r - \alpha } \right)}^p}}}} \\
& \quad+ \sum\limits_{i = 1}^{p - 1} {{{\left( { - 1} \right)}^{i - 1}}\zeta \left( {p + 1 - i,\alpha  + 1} \right)} \sum\limits_{n = 1}^\infty  {\frac{1}{{{{\left( {n + \beta } \right)}^m}{{\left( {n + r - \alpha } \right)}^i}}}} \\
& \quad- \sum\limits_{i = 1}^{m - 1} {{{\left( { - 1} \right)}^{i - 1}}\zeta \left( {m + 1 - i,\beta  + 1} \right)} \sum\limits_{n = 1}^\infty  {\frac{1}{{{{\left( {n + \alpha } \right)}^p}{{\left( {n + r - \beta } \right)}^i}}}} .\tag{3.12}
\end{align*}
Putting $r=2\alpha,\alpha=\beta$ in (3.12), we have that
\begin{align*}
&\left\{ {{{\left( { - 1} \right)}^{m - 1}} - {{\left( { - 1} \right)}^{p - 1}}} \right\}\sum\limits_{n = 1}^\infty  {\frac{{{H_{n + 2\alpha }}}}{{{{\left( {n + \alpha } \right)}^{p + }}^m}}} \\
 =& \left\{ {{{\left( { - 1} \right)}^{m - 1}} -{{\left( { - 1} \right)}^{p - 1}}} \right\}{H_\alpha }\zeta \left( {p + m,\alpha  + 1} \right)\\
 & + \sum\limits_{i = 1}^{p - 1} {{{\left( { - 1} \right)}^{i - 1}}\zeta \left( {p + 1 - i,\alpha  + 1} \right)} \zeta \left( {m + i,\alpha  + 1} \right)\\
 & - \sum\limits_{i = 1}^{m - 1} {{{\left( { - 1} \right)}^{i - 1}}\zeta \left( {m + 1 - i,\beta  + 1} \right)} \zeta \left( {p + i,\alpha  + 1} \right).\tag{3.13}
\end{align*}
From (3.13), we can get some specific cases
\begin{align*}
&\sum\limits_{n = 1}^\infty  {\frac{{{H_{n + 2\alpha }}}}{{{{\left( {n + \alpha } \right)}^3}}}}  = {H_\alpha }\zeta \left( {3,\alpha  + 1} \right) + \frac{1}{2}{\zeta ^2}\left( {2,\alpha  + 1} \right),\\
&\sum\limits_{n = 1}^\infty  {\frac{{{H_{n + 2\alpha }}}}{{{{\left( {n + \alpha } \right)}^5}}}}  = {H_\alpha }\zeta \left( {5,\alpha  + 1} \right) + \zeta \left( {2,\alpha  + 1} \right)\zeta \left( {4,\alpha  + 1} \right) - \frac{1}{2}{\zeta ^2}\left( {3,\alpha  + 1} \right),\\
&\sum\limits_{n = 1}^\infty  {\frac{{{H_{n + 2\alpha }}}}{{{{\left( {n + \alpha } \right)}^7}}}}  = {H_\alpha }\zeta \left( {7,\alpha  + 1} \right) + \zeta \left( {2,\alpha  + 1} \right)\zeta \left( {6,\alpha  + 1} \right) + \frac{1}{2}{\zeta ^2}\left( {4,\alpha  + 1} \right) - \zeta \left( {3,\alpha  + 1} \right)\zeta \left( {5,\alpha  + 1} \right).
\end{align*}
{\bf Acknowledgments.} The authors would like to thank the anonymous
referee for his/her helpful comments, which improve the presentation
of the paper.


{\small
\begin{thebibliography}{99}
\bibitem{A2000} George E. Andrews, Richard Askey, Ranjan Roy. {\sl Special Functions}.
Cambridge University Press., 2000: 481-532.
\bibitem{BBG1994}
David H. Bailey, Jonathan M. Borwein and Roland Girgensohn. {\sl Experimental evaluation of Euler sums}.
Experimental Mathematics., 1994, {\bf 3}(1): 17-30.
\bibitem{B1985} B. C. Berndt. {\sl Ramanujan¡¯s Notebooks, Part I}. Springer-Verlag, New York., 1985.
\bibitem{B1989} B. C. Berndt. {\sl Ramanujan¡¯s Notebooks, Part II}. Springer-Verlag, New York., 1989.
\bibitem{BBG1995}
David Borwein, Jonathan M. Borwein and Roland Girgensohn. {\sl Explicit evaluation of
Euler sums}. Proc. Edinburgh Math., 1995, {\bf 38}: 277-294.
\bibitem{BBGP1996}
J.Borwein, P.Borwein, R.Girgensohn, S.Parnes. {\sl Making sense of experimental mathematics}. Mathematical Intelligencer., 1996, {\bf18}(4): 12-18.
\bibitem{BBBL2001}
Jonathan M. Borwein, David M. Bradley, David J. Broadhurst, Petr. Lison¨§k.
{\sl Special values of multiple polylogarithms. Trans}. Amer. Math. Soc., 2001, {\bf 353}(3): 907-941.
\bibitem{BZB2008}
J. M. Borwein, I. J. Zucker, J. Boersma. {\sl The evaluation of character Euler double sums}. Ramanujan J., 2008, {\bf 15} (3): 377-405.

\bibitem{BG1996}
J.M. Borwein, R. Girgensohn, {\sl Evaluation of triple Euler sums}, Electron. J. Combin., 1996: 2-7.
\bibitem{C2010}
J. Choi, D. Cvijovi$\acute{c}$.{\sl Values of the polygamma functions at rational arguments}. J. Phys. A. Math.
Theor., 2007, {\bf 40}: 15019-15028.
\bibitem{FS1998}
Philippe Flajolet and Bruno Salvy. {\sl Euler sums and contour integral representations}. Experimental Mathematics., 1998, {\bf 7}(1): 15--35.

\bibitem{F2005}
Pedro Freitas. {\sl Integrals of polylogarithmic functions, recurrence relations, and associated Euler sums}. Mathematics of Computation., 2005, {\bf 74}(251): 1425-1440.
\bibitem{H1996}
J. G. Huard, K. S. Williams and N. Y. Zhang. {\sl On Tornheim¡¯s double series}. Acta Arith., 1996, {\bf 75}(2): 105-117.
\bibitem{K1996}
K. K$\ddot{o}$lbig. {\sl The polygamma function $\psi(x)$ for $x = 1/4 $and $x = 3/4$}. J. Comput. Appl. Math., 1996, {\bf 75}: 43-46.
\bibitem{L1974}
Comtet L. {\sl Advanced combinatorics}. D Reidel Publishing Company, Boston., 1974.
\bibitem {M1994}
C. Markett. {\sl Triple sums and the Riemann zeta function}. J. Number Theory., 1994, {\bf 48}(2): 113-132.
\bibitem{M2014}
I. Mez$\ddot{o}$. {\sl Nonlinear Euler sums}. Pacific J. Math., 2014, {\bf 272}: 201-226.
\bibitem {S2003}
A. Sofo. {\sl Computational Techniques for the Summation of Series}. Kluwer Academic/Plenum Pub-
lishers, New York., 2003.
\bibitem {S2009}
A. Sofo. {\sl Integral forms of sums associated with harmonic numbers}. Applied Mathematics and Computation., 2009, {\bf 207}(2): 365-372.
\bibitem{S2010}
 A. Sofo. {\sl Harmonic sums and integral representations}. J. Appl. Anal., 2010, {\bf 16}: 265-277.
\bibitem{S2011}
 A. Sofo, H. M. Srivastava. {\sl Identities for the harmonic numbers and binomial coefficients}. Ramanu-
jan J., 2011, {\bf 25}: 93-113.
\bibitem{So2011}
 A. Sofo. {\sl Harmonic number sums in closed form}. Math. Commun., 2011, {\bf 16}: 335-345.
\bibitem {S2014}
A. Sofo. {\sl Shifted harmonic sums of order two}. Communications of the Korean Mathematical Society., 2014, {\bf 29}(2): 239-255.
\bibitem {S2015}
A. Sofo. {\sl Quadratic alternating harmonic number sums}. J. Number Theory., 2015, {\bf 154}: 144-159.
\bibitem {S2012}
H.M. Srivastava, J. Choi. {\sl Zeta and q-Zeta Functions and Associated Series and Integrals}. Elsevier
Science Publishers., Amsterdam, London and New York, 2012.
\bibitem {S1985}
M. V. Subbarao and R. Sitaramachandra Rao. {\sl On some infinite series of L. J. Mordell and their analogues}. Pacific J. Math., 1985, {\bf 119}(1): 245-255.
\bibitem{Xu2016}
Ce Xu, Jinfa Cheng. {\sl Some Results On Euler Sums}. Functions et Approximatio., 2016, {\bf 54}(1): 25-37.
\bibitem{X2016}
Ce Xu, Yuhuan Yan, Zhijuan Shi. {\sl Euler sums and integrals of polylogarithm functions}. J. Number Theory., 2016, {\bf 165}: 84-108.
\end{thebibliography}}
 \end{document}